\documentclass{fic-l}

\newtheorem{theorem}{Theorem}[section]
\newtheorem{corollary}[theorem]{Corollary}
\newtheorem{lemma}[theorem]{Lemma}
\newtheorem{proposition}[theorem]{Proposition}

\theoremstyle{definition}
\newtheorem{definition}[theorem]{Definition}
\newtheorem{example}[theorem]{Example}
\newtheorem{algorithm}[theorem]{Algorithm}

\theoremstyle{remark}
\newtheorem{remark}[theorem]{Remark}

\numberwithin{equation}{section}

\def\bmath#1{\mbox{\boldmath$#1$}}

\newcommand{\OO}{\bmath{\Omega}}
\newcommand{\DD}{\mathbf{D}}

\newcommand{\CC}{\mathbf{C}}

\newcommand{\PP}{\mathbf{P}}

\newcommand{\TT}{\mathbf{T}}
\newcommand{\UU}{\mathbf{U}}

\newcommand{\XX}{\mathbf{X}}

\newcommand{\YY}{\mathbf{Y}}
\newcommand{\WW}{\mathbf{W}}
\newcommand{\AAA}{\mathbf{A}}

\newcommand{\Kt}{{\widetilde{K}}}

\newcommand{\vecXX}{{\vec{\XX}}}
\newcommand{\vecUU}{{\vec{\UU}}}

\newcommand{\vecYY}{{\vec{\YY}}}
\newcommand{\vecDD}{{\vec{\DD}}}
\newcommand{\vecD}{{\vec{D}}}

\newcommand{\vecU}{{\vec{U}}}
\newcommand{\vecY}{{\vec{Y}}}
\newcommand{\vecX}{{\vec{X}}}

\newcommand{\vecx}{{\vec{x}}}

\newcommand{\vecu}{{\vec{u}}}

\newcommand{\Ac}{\mathcal{A}}
\newcommand{\Bc}{\mathcal{B}}

\newcommand{\Dc}{\mathcal{D}}

\newcommand{\Fc}{\mathcal{F}}

\newcommand{\Lc}{\mathcal{L}}

\newcommand{\Uc}{\mathcal{U}}
\newcommand{\Xc}{\mathcal{X}}

\newcommand{\zh}{\hat{0}}
\newcommand{\oh}{\hat{1}}

\newcommand{\pih}{\hat{\pi}}

\newcommand{\bege}{\begin{equation}}
\newcommand{\ene}{\end{equation}}
\newcommand{\begp}{\begin{proposition}}
\newcommand{\enp}{\end{proposition}}
\newcommand{\begt}{\begin{theorem}}
\newcommand{\ent}{\end{theorem}}
\newcommand{\begl}{\begin{lemma}}
\newcommand{\enl}{\end{lemma}}
\newcommand{\begc}{\begin{corollary}}
\newcommand{\enc}{\end{corollary}}
\newcommand{\begr}{\begin{remark}}
\newcommand{\enr}{\end{remark}}
\newcommand{\begd}{\begin{definition}}
\newcommand{\enf}{\end{definition}}
\newcommand{\begx}{\begin{example}}
\newcommand{\enx}{\end{example}}
\newcommand{\bega}{\begin{array}}
\newcommand{\ena}{\end{array}}

\begin{document}

\title{Extension of Fill's perfect rejection sampling algorithm \\ to
general chains \\
(EXT.\ ABS.)
}

\author{James Allen Fill}
\address{Department of Mathematical Sciences \\
The Johns Hopkins University \\
Baltimore, MD 21218--2682, U.~S.~A.}
\thanks{The first and second authors
have been supported in part by NSF grants DMS--9626756 and
DMS--9803780, and by the Acheson J.~Duncan Fund for the Advancement of Research in
Statistics.}
\email{jimfill@jhu.edu}

\author{Motoya Machida}
\address{Department of Mathematical Sciences \\
The Johns Hopkins University \\
Baltimore, MD 21218--2682, U.~S.~A.}
\email{machida@mts.jhu.edu}

\author{Duncan J.~Murdoch}
\address{Department of Statistical and Actuarial Sciences \\
University of Western Ontario \\
London, Ontario N6G 2E9, Canada} 
\thanks{The third and fourth authors have been supported in part by NSERC}
\email{murdoch@fisher.stats.uwo.ca}

\author{Jeffrey S.~Rosenthal}
\address{Department of Statistics \\
University of Toronto \\
Toronto, Ontario M5S 3G3, Canada}
\email{jeff@math.toronto.edu}

\date{May, 2000}


\begin{abstract}
We provide an extension of the perfect sampling algorithm of Fill (1998) to general
chains, and describe how use of bounding processes can ease computational burden. 
Along the way, we unearth a simple connection between the Coupling From The Past
(CFTP) algorithm originated by Propp and Wilson (1996) and our extension of Fill's
algorithm.
\end{abstract}

\maketitle
 
\section{Introduction}
\label{intro}

Markov chain Monte Carlo (MCMC) methods have become extremely popular
for Bayesian inference problems (consult, e.g.,\ Gelfand and Smith~\cite{GS}, Smith
and Roberts~\cite{SR}, Tierney~\cite{Tierney}, Gilks et al.~\cite{GRS}), and for
problems in other areas, such as spatial statistics, statistical physics, and
computer science (see, e.g.,\ Fill~\cite{Fill} or Propp and Wilson~\cite{PWexact}
for pointers to the literature) as a way of sampling approximately from a
complicated unknown probability distribution~$\pi$.  An MCMC algorithm constructs
a Markov chain with one-step transition kernel~$K$ and stationary
distribution~$\pi$; if the chain is run long enough, then under reasonably weak
conditions (cf.\ Tierney~\cite{Tierney}) it will converge in distribution to~$\pi$,
facilitating approximate sampling.

One difficulty with these methods is that it is difficult to assess convergence to
stationarity.  This necessitates the use of difficult theoretical analysis
(e.g.,\ Meyn and Tweedie~\cite{MT}, Rosenthal~\cite{Rosenthal}) or problematic
convergence diagnostics (Cowles and Carlin~\cite{CC}, Brooks, et al.~\cite{BDR})
to draw reliable samples and do proper inference.

An interesting alternative algorithm, called \emph{coupling from the past\/}
(CFTP), was introduced by Propp and Wilson~\cite{PWexact} (see also~\cite{PWuser}
and~\cite{PWhowto}) and has been studied and used by a number of authors (including
Kendall~\cite{Kendall}, M{\o}ller~\cite{Moller}, Murdoch and Green~\cite{MG}, Foss
and Tweedie~\cite{FT}, Kendall and Th\"onnes~\cite{KT}, Corcoran and
Tweedie~\cite{CT}, Green and Murdoch~\cite{GM}, Murdoch and Rosenthal~\cite{MR}). 
By searching backwards in time until paths from all starting states have coalesced,
this algorithm uses the Markov kernel~$K$ to sample \emph{exactly\/} from~$\pi$.

Another method of perfect simulation, for finite-state stochastically monotone
chains, was proposed by Fill~\cite{Fill}.  Fill's algorithm
is a form of rejection sampling.  This algorithm was later extended by
M{\o}ller and Schladitz~\cite{MS} and Th\"onnes~\cite{Thonnes} to non-finite
chains, motivated by applications to spatial point processes.  Fill's algorithm
has the advantage over CFTP of removing the correlation between the length
of the run and the returned value, which eliminates bias introduced by an impatient
user or a system crash and so is ``interruptible''. However, it has been used only
for stochastically monotone chains, making heavy use of the ordering of state space
elements.  In his paper, Fill~\cite{Fill} indicated that his algorithm could be
suitably modified to allow for the treatment of ``anti-monotone'' chains and (see
his Section 11.2) indeed to generic chains.  In this extended abstract we present a
version of Fill's algorithm for generic chains; we, too, will provide an
explanation in terms of rejection sampling.  We have strived to keep to the spirit
of the talks presented at the Workshop on Monte Carlo Methods held at the Fields
Institute for Research in Mathematical Sciences in Toronto in October, 1998 and
make our results accessible to a broad audience.  Technical details are provided
in the full paper~\cite{FMMR}.

Following is our interruptible algorithm for generic chains.  We discuss some of the
terminology used and other details of the algorithm in Section~\ref{app1}. 

\begin{algorithm}
\label{nontech}
Choose and fix a positive integer~$t$, choose an initial state $\XX_t$ from any
distribution absolutely continuous with respect to~$\pi$, and perform the following
routine.  Run the time-reversed chain~$\Kt$ for~$t$ steps, obtaining $\XX_t,
\XX_{t - 1}, \ldots,$ $\XX_0$ in succession.  Then, reversing the direction of time,
generate (possibly dependent) Markov chains, one [say $\vecYY(x) = (\YY_0(x), \ldots,
\YY_t(x))$, with $\YY_0(x) = x$] evolving from each element~$x$ of the state
space~$\Xc$ and each with transition kernel~$K$.  All these trajectories are to be
multivariately
coupled
\emph{ex post facto\/} with the trajectory $\vecXX = (\XX_0, \ldots, \XX_t)$, which
is regarded as a trajectory from~$K$; in particular, $\vecYY(\XX_0) = \vecXX$. 
Finally, we check whether all the values $\YY_t(x)$, $x \in \Xc$, agree.
If they do, we call this
\emph{coalescence\/}, and the value~$\XX_0$ is accepted as an observation
from~$\pi$.  If not, then the value~$\XX_0$ is rejected and so the routine
fails.  We then start the routine again with an independent simulation,
perhaps
with a fresh choice of~$t$ and~$\XX_t$,
and repeat until the algorithm succeeds.
\end{algorithm}

Here is a
simple intuitive
explanation of why the algorithm works correctly.
Imagine (1)~\emph{starting\/} the
construction with $\XX_0 \sim \pi$ and independently (2)~building all of the paths
$\vecYY(x)$ [including $\vecXX = \vecYY(\XX_0)$] \emph{simultaneously\/}.  Since
determination of coalescence and the value of the coalesced paths at
time~$t$ rely only on the second piece of randomness, conditionally given
coalescence to $\XX_t = z$ (for any~$z$) we will still have
$\XX_0 \sim \pi$, as desired.  The algorithm builds the randomness in a different
order, starting from~$\XX_t$.
%

\begin{remark}
(a)~Note that no assumption is made in Algorithm~\ref{nontech} concerning
monotonicity or discreteness of the state space.

(b)~See~(\ref{reversal}) for the definition of~$\Kt$.

(c)~To couple all of the trajectories~$\vecYY(x)$ \emph{ex post facto\/} with
the trajectory $\vecXX$ means first to devise a multivariate coupling for
all the trajectories by means of a transition rule and then to employ that coupling
conditionally having observed the single trajectory $\vecYY(\XX_0) = \vecXX$.  Just
how this is done is described in Section~3.

(d)~If coalescence occurs, then of course the common value of~$\YY_t(x)$, $x \in \Xc$,
is the initial state~$\XX_t$.

(e)~We have reversed the direction of time, and the roles of the kernels~$K$
and~$\Kt$, compared to Fill~\cite{Fill}.  Furthermore, Fill's original algorithm
also incorporated a search for a good value of~$t$ by doubling the previous value
of~$t$ until the first success.  For the most part, we shall not address such
issues, instead leaving the choice of~$t$ entirely up to the user; but see
Section~\ref{altalgsub}.

(f)~We have included the technical absolute continuity restriction on the
distribution of~$\XX_t$ to ensure correctness.  In a typical application, one might
know a density for~$\pi$ with respect to some measure~$\lambda$ (for example,
$\lambda =$ Lebesgue measure) up to a normalizing constant.  Then if, for example,
that density is positive on the entire state space~$\Xc$, one need only take~$\XX_t$
to have any distribution having density with respect to~$\lambda$ (for example, a normal
distribution). 

Also, if the state space is discrete, or if it is continuous and all probabilistic
ingredients to the algorithm (such as the kernel~$K$) are sufficiently smooth,
then the user may choose~$\XX_t$ deterministically and arbitrarily.
\end{remark}

As mentioned above, we
will discuss the details of Algorithm~\ref{nontech} in Section~\ref{app1}.  First,
in Section~\ref{motivation}, we motivate our algorithm in the context of a rather
general rejection sampling framework.  A more rigorous treatment may be found in the 
full paper~\cite{FMMR}.

In Section~\ref{app2} we discuss how the computational burden of tracking all of
the trajectories~$\YY(x)$ can be eased by the use of coalescence detection events
in general and bounding processes in particular; these processes take on a very
simple form (see Section~\ref{mono}) when the state space is partially
ordered and the transition rule employed is monotone.  In the full paper we present
a computationally efficient modification of Algorithm~\ref{nontech} that applies
when~$K$ is assumed to be stochastically monotone.  (As discussed in Fill and
Machida~\cite{SMRM} and Machida~\cite{Machida}, this is a weaker assumption than
that of the existence of a monotone transition rule.)  When the state space is
finite, the algorithm for the stochastically monotone case reduces to Fill's
original algorithm.  The full paper also discusses a ``cross-monotone''
generalization of stochastic monotonicity.
%

In Section~\ref{taleof2} we compare Algorithm~\ref{nontech} and CFTP.  We also
demonstrate
a simple connection between CFTP and an infinite-time-window version of
Algorithm~\ref{nontech} (namely, Algorithm~\ref{altalg}).

This extended abstract discusses only algorithms for generating a single
observation from a distribution~$\pi$ of interest.  In the full paper, we discuss
various strategies for perfect generation of samples of
arbitrary size from~$\pi$. 

The goal of this extended abstract is to describe how to apply Fill's perfect
sampling algorithm to a broad spectrum of problems of real interest, not to
develop any specific applications in detail.  But an underlying theme here is
that while our extended algorithm (Algorithm~\ref{nontech}) tends to be computationally
more intricate than CFTP (see Section~\ref{comparison}), theoretically it is as broadly
applicable as is CFTP.  In this spirit, we will point to the applied perfect sampling
literature at appropriate junctures, taking note of past applications of both CFTP and
Fill's algorithm as examples.  We hope that our extension of the latter algorithm will
stimulate further research into this less-used alternative for perfect MCMC
simulation.

A valuable background resource is the annotated bibliography of perfect sampling
maintained by Wilson~\cite{Wilson}. 

%
%

\section{Rejection sampling using auxiliary \mbox{randomness}}
\label{motivation}

Given a bivariate distribution $\Lc(X, X')$, suppose
that, for each~$x'$, we can simulate~$X$ from the conditional distribution~$\Lc(X
| X' = x')$.  How can we simulate~$X$ from its marginal distribution $\pi :=
\Lc(X)$?  This is the problem confronting us in the context of
Algorithm~\ref{nontech}, with $(X, X') := (\XX_0, \XX_t)$ and $\XX_t$ chosen
according to~$\pi$.  Indeed, we assumed there that the user can simulate from
$\Lc(\XX_0\,|\,\XX_t = x') = \Kt^t(x', \cdot)$, where $\Kt^t$ is the $t$-step
time-reversal transition kernel; and, when $\XX_t \sim \pi$, we have $\Lc(\XX_0)
= \pi$.  Of course, if the user could simulate~$\XX_t \sim \pi$, then
either~$\XX_t$ or~$\XX_0$ could be used directly as an observation from~$\pi$. 
But, for MCMC, this is an unreasonable assumption.

So we turn to rejection sampling (e.g.,\ Devroye~\cite{Devroye}),
done conditionally given $X' = x'$.  Our goal is to use the observation $x$ from
$\Lc(X | X' = x')$~to simulate one from~$\pi$.  This can be done by
accepting~$x$ as an observation from~$\pi$ with probability
$$
\gamma_{x', x} := \alpha_{x'} \frac{\pi(dx)}{P(X \in dx | X' = x')}
$$
for any $\alpha_{x'}$ chosen to make $0 < \gamma_{x', x} \leq 1$;
indeed, then
$$
P(X \in dx\,|\,X' = x', \mbox{accept}) = \frac{P(X \in dx, \mbox{accept}\,|\,X'
= x')}{P(\mbox{accept}\,|\,X' = x')} = \frac{\alpha_{x'} \pi(dx)}{\alpha_{x'}
\int\!\pi(dy)} = \pi(dx).
$$
The question remains how to engineer a coin-flip with probability~$\gamma_{x',x}$
of heads, given that one can rarely \emph{compute\/}~$\gamma_{x', x}$ in practice.

However, if we can find an event~$C$ so that
\begin{equation}
\label{ceq}
P(C\,|\,X' = x', X = x) = \gamma_{x', x}\mbox{\ \ \ for all $x', x$,}
\end{equation}
then we need only accept when (and only when) $C$ occurs.  Condition~(\ref{ceq})
requires precisely that
\begin{equation}
\label{eq2}
\ \ \ P(\{X' \in B'\} \cap C \cap \{X \in B\}) = \pi(B) \int_{B'}\,\alpha_{x'}
P(X' \in dx')\mbox{\ \ \ for all $B, B'$.}
\end{equation}
Note that if we can choose~$C$ so that condition~(\ref{eq2}) holds, then
setting~$B$ to be the entire state space~$\Xc$ for~$X$, we obtain
$$
P(\{X' \in B'\} \cap C) = \int_{B'}\,\alpha_{x'} P(X' \in dx')\mbox{\ \ \
for all $B'$}
$$
and hence
\begin{equation}
\label{fact}
\ \ \ P(\{X' \in B'\} \cap C \cap \{X \in B\}) = \pi(B)\,P(\{X' \in B'\} \cap
C) \mbox{\ \ \ for all $B, B'$.}
\end{equation}
Conversely, if we can choose~$C$ so that condition~(\ref{fact}) holds, then we
can choose
$$
\alpha_{x'} := P(C\,|\,X' = x')
$$
to satisfy~(\ref{ceq}) and~(\ref{eq2}).

How might we choose~$C$ to satisfy~(\ref{fact})?  Observe that if~$C$ and another
variable~$Y$ are such that (i)~$Y = X'$ over the event~$C$ and (ii)~$X$ and the
pair~$(Y, C)$ are independent, then
\begin{eqnarray*}
\mbox{LHS(\ref{fact})}
 &=& P(\{X' \in B'\} \cap C \cap \{X \in B\}) = P(\{Y \in B'\} \cap C \cap \{X
       \in B\}) \\
 &=& P(\{Y \in B'\} \cap C)\,\pi(B).
\end{eqnarray*}
Setting $B = \Xc$ and then substituting, we obtain~(\ref{fact}).  The
Algorithm~\ref{nontech} is designed precisely so that we may take $(Y, C) :=
(\YY_t(x^*_0), \{\mbox{coalescence}\})$ for an arbitrarily chosen (but fixed)
$x^*_0 \in \Xc$ to satisfy~(i) and~(ii).
This is discussed further in Section~\ref{rigor} below.

\section{Details for Algorithm~1.1}
\label{app1}

\subsection{The ingredients}
\label{rigor}

The goal of this subsection 
is to describe in some detail how to apply
Algorithm~\ref{nontech}.
\medskip

%
\emph{The space $(\Xc, \Bc)$:\/}\ It is sufficient that $(\Xc, \Bc)$ be a complete
separable metric space.  In particular, this covers the case that~$\Xc$ is discrete or
Euclidean.  See Section~3.1 of the full paper~\cite{FMMR} for further discussion. 
\medskip

\emph{The kernel~$K$ and its time-reversal~$\Kt$:\/}\ Let~$K$ be a Markov transition
kernel on~$\Xc$.  The kernel is chosen (by the user) so that~$\pi$ is a stationary
distribution, i.e.,\ so that
$$
\int_{\Xc} \pi(dx) K(x, dy) = \pi(dy)\mbox{\ \ \ on $\Xc$.}
$$
The time-reversed kernel~$\Kt$ (also on~$\Xc$) is then defined by
\begin{equation}
\label{reversal}
\pi(dx) K(x, dy) = \pi(dy) \Kt(y, dx)\mbox{\ \ \ on $\Xc \times \Xc$.}
\end{equation}

\emph{The transition rule~$\phi$:\/}\ There exists a transition
rule which can be used to drive the construction of the Markov chain of interest. 
More precisely, there exists a probability space $(\Uc, \Fc, \mu)$
and a (suitably measurable) function
$\phi: \Xc \times \Uc \to \Xc$
such that
\begin{equation}
\label{kq}
K(x, B) = \mu\{u: \phi(x, u) \in B\},\ \ \ x \in \Xc,\ \ B \in \Bc.
\end{equation}
Such~$\phi$ (with accompanying~$\mu$) is sometimes called a \emph{transition rule\/}.
We choose and fix such a $(\phi, \mu)$.

\begin{remark}
\label{phiremark}
%
(a)~A transition rule~$\phi$ can always be found that uses $(\Uc, \Fc, \mu) =
([0, 1]$, $\mbox{Borels}, \mbox{uniform distribution})$.
In the special case that~$\Xc$ is the unit interval,
we can in fact use
$$
\phi(x, u) \equiv G(x, u)
$$
where~$G(x, \cdot)$ is the usual inverse probability transform corresponding to
the distribution function $u \mapsto K(x, [0,u])$.

(b)~If~$\Xc$ is discrete (finite or countably infinite),
a very simple alternative choice is the following
``independent-transitions'' transition rule.  Let~$\Uc = \Xc^{\Xc}$,
let~$\mu$ be product measure with $x$th marginal $K(x, \cdot)$ ($x \in \Xc$), and
let~$\phi$ be the evaluation function
$$
\phi(x, u) := u(x).
$$

(c) Many interesting examples of transition rules can be found in the literature,
including Diaconis and Freedman~\cite{DF} and the references cited in
Section~\ref{intro}.

(d) Usually there is a wealth of choices of transition rule, and the art is to
find one giving rapid and easily detected coalescence.  Without going into details
at this point, we remark that the transition rule in~(b) usually performs quite
badly, while transition rules having a certain monotonicity property will perform
well under monotonicity assumptions on~$K$.
\end{remark}

\emph{The Markov chain and a first probability space:\/}\ 
We will need to set up two probability spaces.
The first space $(\Omega, \Ac, P)$---designated in
ordinary typeface---will be useful for theoretical considerations and for the
computation of certain conditional probability distributions.
The second
space $(\OO, \AAA, \PP)$---designated in boldface type---will be the probability
space actually simulated when the algorithm is run.
All random variables defined on the first space (respectively, second space) will
also be designated in ordinary typeface (resp.,\ boldface type).  We have chosen
this notational system to aid the reader:\ Corresponding variables, such as~$X_0$
and~$\XX_0$, will play analogous roles in the two spaces.

From our previous
comments it is easy to see that there exists a
space~$(\Uc, \Fc)$, a transition rule $(\phi, \mu)$, and a probability space~$(\Omega,
\Ac, P)$ on which are defined independent random variables $X_0, U_1, U_2, \ldots,
U_t$ with $X_0 \sim \pi$ and each $U_s \sim \mu$.  Now inductively define
\begin{equation}
\label{drive}
X_s := \phi(X_{s - 1}, U_s),\ \ \ 1 \leq s \leq t.
\end{equation}
Then $\vecX := (X_0, \ldots, X_t)$ is easily seen to be a stationary Markov
chain with kernel~$K$, in the sense that 
\begin{equation}
\label{ftraj}
\ \ \ \ \ \ \ P(X_0 \in dx_0, \ldots, X_t \in dx_t) = \pi(dx_0) K(x_0,
dx_1)\!\cdots\!K(x_{t-1}, dx_t)\mbox{\ \ \ on~$\Xc^{t+1}$.}
\end{equation}
In fact, for each $x \in \Xc$ we obtain a chain with kernel~$K$ started from~$x$
by defining $Y_0(x) := x$ and, inductively,
$$
Y_s(x) := \phi(Y_{s-1}(x), U_s).
$$
Let $\vecY(x) := (Y_0(x), \ldots, Y_t(x))$.  In this notation we have
$\vecY(X_0) = \vecX$.  Let
\begin{equation}
\label{Cdef}
C := \{Y_t(x)\mbox{\ does not depend on~$x$}\}
\end{equation}
denote the
event that the trajectories~$\vecY(x)$
have all coalesced by time~$t$.
Fixing $x^*_0 \in \Xc$ arbitrarily and taking
$$
X = X_0,\ \ \ X' = X_t,\ \ \ Y = Y_t(x^*_0),\ \ \ \mbox{$C$ as at~(\ref{Cdef})}
$$
to match up with the notation in Section~\ref{motivation}, key observations are that,
for any (measurable) subset~$B$ of~$\Xc$, $Y = X'$ over the event~$C$, and~$X$ and the
event $\{Y \in B\} \cap C$ are independent.  The independence here follows from the fact
that~$X_0$ and $\vecU := (U_1, \ldots, U_t)$ have been chosen to be independent.
%
\medskip

\emph{A second probability space and the algorithm:\/}\ We make use of
the auxiliary randomness provided by $X_1,
\ldots, X_{t - 1}$ and $\vecU$ and compute conditional probability distributions
for the first probability space in stages.  First observe
from~(\ref{ftraj}) and repeated use of~(\ref{reversal}) that
$$
P(X_0 \in dx_0, \ldots, X_{t - 1} \in dx_{t - 1}\,|\,X_t = x_t) = \Kt(x_t, dx_{t -
1}) \cdots \Kt(x_1, dx_0).
$$
Next, we will discuss in Section~\ref{imp} how to compute
$\Lc(\vecU\,|\,\vecX = \vecx)$.  Finally,
whether or not~$C$ occurs is determined solely by the randomness in~$\vecU$, so the
conditional probability of~$C$ given $(\vecX, \vecU) = (\vecx, \vecu)$ is degenerately
either~$0$ or~$1$.  
 
Moreover, our discussion has indicated how to set up and \emph{simulate\/} the
second space.  As discussed in Section~\ref{intro},
we assume
that the user knows enough about~$\pi$ qualitatively to be able to
simulate~$\XX_t$ so that $\Lc_{\PP}(\XX_t) \ll \pi$.
Having chosen $\XX_t = x_t$, the user draws an observation $\XX_{t - 1} = x_{t - 1}$ from
$\Kt(x_t, \cdot)$, then an observation $\XX_{t - 2} = x_{t - 2}$ from $\Kt(x_{t - 1},
\cdot)$, etc.\ \ Next, having chosen $\vecXX = \vecx$ [i.e., $(\XX_0, \ldots,
\XX_t) = (x_0, \ldots, x_t)$], the user draws an observation $\vecUU = \vecu$
from $\Lc(\vecU\,|\,\vecX = \vecx)$
and constructs $\vecYY(x) = (\YY_0(x), \ldots, \YY_t(x))$ by setting $\YY_0(x) := x$
and, inductively,
$$
\YY_s(x) := \phi(\YY_{s - 1}(x), \UU_s).
$$
Finally, the user
declares that~$\CC$, or \emph{coalescence\/}, has occurred if and only if
the values of~$\YY_t(x)$, $x \in \Xc$, all agree.
The
conditional distribution
of output from Algorithm~\ref{nontech} given that it ultimately succeeds (perhaps
only after many iterations of the basic routine) is~$\pi$, as desired.

\begin{remark}
\label{genremark}
(a)~If $\PP(\CC) > 0$ for suitably large~$t$, then ultimate success is (a.s.)\
guaranteed if the successive choices of~$t$ become large.  A necessary condition
for ultimate positivity of~$\PP(\CC)$ is uniform ergodicity of~$K$.  This condition
is also sufficient, in the (rather weak) sense that if~$K$ is uniformly ergodic,
then there exists a finite integer~$m$ and a transition rule~$\phi_m$ for
the $m$-step kernel~$K^m$ such that Algorithm~\ref{nontech}, applied
using~$\phi_m$, has $\PP(\CC) > 0$ when~$t$ is chosen sufficiently large.
Compare the analogous Theorem~4.2 for CFTP in Foss and Tweedie~\cite{FT}.

(b)~Just as discussed in Fill~\cite{Fill} (see especially the end of Section~7
there), the algorithm (including its repetition of the basic routine) we have
described is interruptible; that is, its running time (as measured by number of
Markov chain steps) and output are independent random variables, conditionally
given that the algorithm eventually terminates.

(c)~If the user chooses the value of~$\XX_t$ ($= z$, say) deterministically, then
all that can be said in general is that the algorithm works properly for
$\pi$-a.e.\ such choice.  In this case, let the notation $\PP_z(\CC)$ reflect the
dependence of~$\PP(\CC)$ on the initial state~$z$.  Then clearly
$$
\int\,\PP_z(\CC)\,\pi(dz) = P(C),
$$ 
which is the unconditional probability of coalescence in our first probability
space and therefore equal to the probability that CFTP terminates over an interval
of width~$t$.  This provides a first link between CFTP and
Algorithm~\ref{nontech}.  (Very) roughly recast, the distribution of running time
for CFTP is the stationary mixture, over initial states, of the distributions of
running time for Algorithm~\ref{nontech}.  For further elaboration of the
connection between the two algorithms, see Section~\ref{conn}.
\end{remark}

\subsection{Imputation}
\label{imp}

In order to be able to run Algorithm~\ref{nontech}, the user needs to be able to
impute~$\vecU$ from~$\vecX$, i.e.,\ to draw from $\Lc(\vecU\,|\,\vecX = \vecx)$. 
In this subsection we explain how to do this.  We begin with the computation
%
\begin{eqnarray*}
\lefteqn{P(\vecU \in d\vecu\,|\,\vecX = \vecx)} \\
 & =& P(\vecU \in d\vecu\,|\,X_0 = x_0,\,\phi(x_0, U_1) = x_1,\,\ldots,\,\phi(x_{t
        - 1}, U_t) = x_t)\mbox{\ \ \ by~(\ref{drive})} \\
 & =& P(\vecU \in d\vecu\,|\,\phi(x_0, U_1) = x_1,\,\ldots,\,\phi(x_{t - 1}, U_t) =
        x_t)\mbox{\ \ \ by indep.\ of~$X_0$ and~$\vecU$} \\
 & =& P(U_1 \in du_1\,|\,\phi(x_0, U_1) = x_1) \times \cdots \times P(U_t \in
        du_t\,|\,\phi(x_{t - 1}, U_t) = x_t) \\
 &{}& \mbox{\ \ \ \ \ \ \ by independence of $U_1, \ldots, U_t$} \\
 & =& P(U_1 \in du_1\,|\,\phi(x_0, U_1) = x_1) \times \cdots \times P(U_1 \in
        du_t\,|\,\phi(x_{t - 1}, U_1) = x_t) \\
 &{}& \mbox{\ \ \ \ \ \ \ since $U_1, \ldots, U_t$ are identically distributed} \\
 & =& P(U_1 \in du_1\,|\,X_0 = x_0,\,X_1 = x_1) \times \cdots \times P(U_1 \in
        du_t\,|\,X_0 = x_{t - 1},\,X_1 = x_t),
\end{eqnarray*}
where the last equality is justified in the same fashion as for the first two.
%
This establishes
\begin{lemma}
\label{condlemma}
The $t$-fold product of the measures
$$
P(U_1 \in du_1\,|\,X_0 = x_0,\,X_1 = x_1),\ \ldots,\  P(U_1 \in du_t\,|\,X_0 = x_{t
  - 1},\,X_1 = x_t)
$$
serves as a conditional probability distribution $P(\vecU \in d\vecu\,|\,\vecX =
\vecx)$.
\end{lemma}

In setting up the second probability space, therefore, the user, having chosen
$\vecXX = \vecx$, draws an observation $\vecUU = \vecu$ by drawing $\UU_1, \ldots,
\UU_t$ independently, with $\UU_s$ chosen according to the distribution
$\Lc(U_1\,|\,X_0 = x_{s - 1},\,X_1 = x_s)$.

\begin{remark}
\label{subtle}
%
(a)~If~$\Xc$ is discrete, suppose we use the ``independent-transitions'' rule~$\phi$
discussed
in Remark~\ref{phiremark}(b).  Then the measure~$\mu$, but with the $x_0$th marginal
replaced by~$\delta_{x_1}$, serves as $\Lc(U_1\,|\,\phi(x_0, U_1) = x_1) =
\Lc(U_1\,|\,U_1(x_{0}) = x_1)$ and therefore as $\Lc(U_1\,|\,X_0 = x_0,\,X_1
= x_1)$.  Informally stated, having chosen $\XX_s = x_s$ and $\XX_{s - 1} = x_{s -
1}$, the user imputes the forward-trajectory transitions from time~$s - 1$ to
time~$s$ in Algorithm~\ref{nontech} by declaring that the transition from
state~$x_{s - 1}$ is to state~$x_s$ and that the transitions from other states are
chosen independently according to their usual non-$\vecX$-conditioned
distributions.

(b)~As another example, suppose that $\Xc = [0, 1]$ and we use the inverse
probability transform transition rule discussed in Remark~\ref{phiremark}(a).
Suppose also that each distribution function $F(x_0, \cdot) = K(x_0, [0, \cdot])$
is strictly increasing and onto~$[0, 1]$.
Then~$\delta_{F(x_0, x_1)}$ serves as $\Lc(U_1\,|\,X_0 = x_0,\,X_1 = x_1)$.  Informally
stated, a generated pair $(\XX_s, \XX_{s - 1}) = (x_s, x_{s-1})$ completely determines
the value $F(x_{s - 1}, x_s)$ for~$\UU_s$.
\end{remark}

\subsection{A toy example}
\label{toy}

We illustrate Algorithm~\ref{nontech} for a very simple example and two different
choices of transition rule.  Consider the discrete state space
$ \Xc = \{0, 1, 2\}$,
and let~$\pi$ be uniform on~$\Xc$.  Let~$K$ correspond to
simple symmetric random walk with holding probability~$1/2$ at the endpoints; that
is, putting $k(x, y) := K(x, \{y\})$,
\begin{eqnarray*}
&& k(0, 0) = k(0, 1) = k(1, 0) = k(1, 2) = k(2, 1) = k(2, 2) = 1/2, \\
&& k(0, 2) = k(1, 1) = k(2, 0) = 0.
\end{eqnarray*}
The stationary distribution is~$\pi$.  As for any ergodic birth-and-death chain,
$K$ is reversible with respect to~$\pi$, i.e.,\ $\Kt = K$.  Before starting the
algorithm, choose a transition rule; this is discussed further below.

For utter simplicity of description, we choose~$t = 2$ and (deterministically)
$\XX_t = 0$ (say); as discussed
in Section~\ref{intro}, a
deterministic start is permissible here.  We then choose $\XX_1 \sim K(0, \cdot)$
and $\XX_0\,|\,\XX_1 \sim K(\XX_1, \cdot)$.  How we proceed from this
juncture depends on what we chose for~$\phi$.

One choice is the independent-transitions rule discussed in
Remarks~\ref{phiremark}(b) and~\ref{subtle}(a).  The algorithm's routine can then
be run using~$6$ independent random bits: these decide~$\XX_1$ (given~$\XX_2$),
$\XX_0$ (given~$\XX_1$), and the~$4$ transitions in the second (forward)
phase of the routine not already determined from the rule
$$
\XX_{s - 1} \mapsto \XX_s \mbox{\ from time $s - 1$ to time~$s$ ($s = 1, 2$).}
$$
There are thus a total of~$2^6 = 64$ possible
overall simulation results, each having probability~$1/64$.  We check that
exactly~$12$ of these produce coalescence.  Of these~$12$ accepted results,
exactly~$4$ have $\XX_0 = 0$, another~$4$ have $\XX_0 = 1$, and a final~$4$
have $\XX_0 = 2$.  Thus $\PP(\CC) = 12/64 = 3/16$, and we confirm that
$\Lc_{\PP_{\CC}}(\XX_0) = \pi$,
as should be true.  An
identical result holds if instead we choose $\XX_t = 1$ or $\XX_t = 2$.

An alternative choice adapts Remarks~\ref{phiremark}(a) and~\ref{subtle}(b) to our
discrete setting.  Note that we can use
$(\Uc, \mu) = (\{0, 1\}, \mbox{uniform})$
and
$$
\phi(\cdot, u) = \left\{
 \begin{array}{lll}
  \mbox{the mapping taking\ \ $0,1,2$ to} & \mbox{$0,0,1$,} & \mbox{respectively,
                                                                   if $u = 0$}
\vspace{.1in} \\
  \mbox{the mapping taking\ \ $0,1,2$ to} & \mbox{$1,2,2$,} & \mbox{respectively,
                                                                   if $u = 1$.}
 \end{array}
 \right.
$$
Choosing $t = 2$ and $\XX_t = 0$ as before, the algorithm can now be run with
just~$2$ random bits.  In this case we check that exactly~$3$ of the~$4$ possible
simulation results produce coalescence, $1$ each yielding $\XX_0 = 0,1,2$.  Note
that $\PP(\CC) = 3/4$ is much larger for this choice of~$\phi$.  In fact,
since~$\phi$ is a monotone transition rule (see Definition~4.2 in Fill~\cite{Fill}
or Definition~\ref{mondef} below), for the choice $\XX_t = 0$ it gives the highest
possible value of~$\PP(\CC)$ among all choices of~$\phi$: see Remark~9.3(e) in
Fill~\cite{Fill}.  It also is a best choice when $\XX_t = 2$.  [On a minor negative
note, we observe that $\PP(\CC) = 0$ for the choice $\XX_t = 1$.  Also note that
the $\pi = (1/3, 1/3, 1/3)$-average of the acceptance probabilities~$(3/4, 0,
3/4)$, namely, $1/2$, is the probability that forward coupling (or CFTP) done with
the same transition rule gives coalescence within~$2$ time units; this corroborates
Remark~\ref{genremark}(c).]

\begin{remark}
\label{toyrate}
Both choices of~$\phi$ are easily extended to handle simple symmetric random
walk on $\{0, \ldots, n\}$ for any~$n$; the second (monotone) choice is again
best possible.  (We assume $\XX_t = 0$.)  For fixed $c \in (0, \infty)$
and large~$n$, results in Fill~\cite{Fill} and Section~4 of Diaconis and
Fill~\cite{DFAP} imply that, for $t = c n^2$, the routine's success probability is
approximately~$p(c)$; here~$p(c)$ increases smoothly from~$0$ to~$1$ as $c$
increases from~$0$ to~$\infty$.  We have not attempted the corresponding
asymptotic analysis for the independent-transitions rule.  Of course our chain is
only a ``toy'' example anyway, because direct sampling from~$\pi$ is elementary.
\end{remark}

\section{Coalescence detection and bounding processes}
\label{app2}

\subsection{Conservative detection of coalescence and detection processes}
\label{detection}

Even for large finite state spaces~$\Xc$, determining exactly whether or not
coalescence occurs in Algorithm~\ref{nontech} can be prohibitively expensive
computationally; indeed, in principle this requires tracking each of the
trajectories $\vecY(x)$, $x \in \Xc$ to completion.  However, observe that if we
repeat the development of Sections~\ref{rigor}--\ref{imp}, replacing the
coalescence event $\{\mbox{$Y_t(x)$ does not depend on~$x$}\}$ of~(\ref{Cdef}) by
any
\emph{subset\/} $C$ of this event whose occurrence (or not) can still be determined
solely from~$\vecU$,
then everything goes
through as before.  We call such an event~$C$ a \emph{coalescence detection
event\/} and reiterate that~$C$ is a conservative indication of coalescence:\ the
occurrence of a given coalescence detection event is sufficient, but not
necessary, for the occurrence of coalescence of the paths
$\vecY(x)$.

In practice, a coalescence detection event is constructed in terms of
a \emph{detection process\/}.  What we mean by this is a stochastic process $\vecD =
(D_0, \ldots, D_t)$, defined on the same probability space $(\Omega, \Ac, P)$
as~$\vecU$ and~$\vecX$, together with a
subset~$\Delta$ of its state space~$\Dc$, such that
\begin{enumerate}
\item[(a)] $\vecD$ is constructed from~$\vecU$, and
\item[(b)] $\{D_s \in \Delta\mbox{\ for some $s \leq t$}\} \subseteq
\{\mbox{$Y_t(x)$ does not depend on~$x$}\}$.
\end{enumerate}
Then $C := \{D_s \in \Delta\mbox{\ for some $s \leq t$}\}$ is a coalescence
detection event.

\begin{remark}
\label{detectionremark}
In practice,
$\vecD$ usually evolves Markovianly using~$\vecU$; more precisely, it is
typically the case that there exists deterministic $d_0 \in \Dc$ and
$\delta: \Dc \times \Uc \to \Dc$ such that $D_0 = d_0$ and
[paralleling~(\ref{drive})]
$$
D_s := \delta(D_{s-1}, U_s),\ \ \ 1 \leq s \leq t.
$$
\end{remark} 

The important consequence is that, having determined the trajectory~$\vecXX$ and
the imputed~$\vecUU$, the user need only follow a single trajectory in the forward
phase of the routine, namely, that of~$\vecD$ (or rather its analogue~$\vecDD$ in
the simulated probability space).

\begin{example}
\label{MTFandtrees}
We sketch two illustrative examples of the use of detection processes that do
not immediately fall into the more specific settings of Sections~\ref{bounding}
or~\ref{mono}.  We hasten to point out, however, that because of the highly special
structure of these two examples, efficient implementation of
Algorithm~\ref{nontech} avoids the use of the forward phase altogether; this is
discussed for example~(a) in Fill~\cite{MTFalgo}.

(a)~Our first example is provided by the move-to-front (MTF) rule studied
in~\cite{MTFalgo}.  Let~$K$ be the Markov kernel corresponding to~MTF with
independent and identically distributed record requests corresponding to
probability weight vector $(w_1, \ldots, w_n)$; see~(2.1) of~\cite{MTFalgo} for
specifics.  The arguments of Section~4 of~\cite{MTFalgo} show that if~$D_s$ is
taken to be the set of all records requested at least once among the first~$s$
requests and~$\Delta$ is taken to consist of all $(n - 1)$-element subsets of the
records $1, \ldots, n$, then~$\vecD$ is a detection process.  Similar detection
processes can be built for the following generalizations of MTF:\ move-to-root for
binary search trees (see Dobrow and Fill~\cite{DF1}~\cite{DF2}) and MTF-like
shuffles of hyperplane arrangement chambers and more general structures (see
Bidigare, et al.~\cite{BHR} and Brown and Diaconis~\cite{BD}).

(b)~A second example of quite similar spirit is provided by the (now well-known)
Markov chain~$(X_t)$ for generating a random spanning arborescence of the
underlying weighted directed graph, with vertex set~$\Uc$, of a Markov chain
$(U_t)$ with kernel~$q$.  Consult Propp and Wilson~\cite{PWhowto} (who also
discuss a more efficient ``cycle-popping'' algorithm) for details.  We consider
here only the special case that~$(U_t)$ is an i.i.d.\ sequence, i.e.,\ that $q(v,
w) \equiv q(w)$.  A transition rule~$\phi$ for the chain~$(X_t)$ is created as
follows:\ for vertex~$u$ and arborescence~$x$ with root~$r$, $\phi(x, u)$ is the
arborescence obtained from~$x$ by adding an arc from~$r$ to~$u$ and deleting the
unique arc in~$x$ whose tail is~$u$.  Then it can be shown that if $D_s$ is taken
to be the set of all vertices appearing at least once in $(U_1, \ldots, U_s)$ and
$\Delta := \{\Uc\}$, then~$\vecD$ is a detection process.
\end{example}

\subsection{Bounding processes}
\label{bounding}

We obtain a natural example of a detection process~$\vecD$ when
(a)~$\vecD$ is constructed from~$\vecU$,
(b)~the corresponding state space~$\Dc$ is some
collection of subsets of~$\Xc$, with
$$
\Delta := \{ \{x'\}:\ x' \in \Xc\},
$$
and
$$
\mbox{(c)\ \ \ }D_s \supseteq \{Y_s(x):\ x \in \Xc\}.
$$
The concept is simple: in this case, each set~$D_s$ is just a ``conservative
estimate'' (i.e., a superset) of the corresponding set $\{Y_s(x): x \in \Xc\}$
of trajectory values; thus if $D_s = \{x'\}$, then the trajectories~$\vecY(x)$
are coalesced to state~$x'$ at time~$s$ and remain coalesced thereafter.  We
follow the natural impulse to call such a set-valued detection process
a \emph{bounding process\/}.  Such bounding processes arise naturally in the contexts of
monotone and anti-monotone transition rules (and have been used by many authors):\
see the next subsection.  Other examples of bounding processes can be found in
works of Huber:\ see~\cite{Hefficient} and~\cite{Hexact} in connection with CFTP
and~\cite{Hinterruptible} in connection with our algorithm.

Of course, nothing is gained, in comparison to tracking all the trajectories, by
the use of a bounding process unless the states of~$\Dc$ have more concise
representations than those of generic subsets of~$\Xc$; after all,
we could always choose $\Dc = 2^{\Xc}$ and $D_s =
\{Y_s(x): x \in \Xc\}$.  One rather general setting where compact representations
are often possible, discussed in the next subsection, is that of a partially
ordered set (poset)~$\Xc$.

\subsection{Monotone transition rules}
\label{mono}

We now suppose that
$\Xc$ is equipped with a partial order.
We also assume here that there exist (necessarily unique) elements~$\zh$ and~$\oh$
in~$\Xc$ (called \emph{bottom element\/} and \emph{top element\/}, respectively) such that
$\zh \leq x \leq \oh$ for all $x \in \Xc$.  We will discuss the case of monotone
transition rules, where we can build
from~$\vecU$
a bivariate
process
$((L_s, V_s))$, taking values at each time~$s$ in $\Xc \times \Xc$, such that
\begin{equation}
\label{squeeze}
L_s \leq Y_s(x) \leq V_s\mbox{\ \ \ for all $0 \leq s \leq t$ and all $x \in
\Xc$.}
\end{equation}
Then $D_s := [L_s, V_s] = \{x \in \Xc:\ L_s \leq x \leq V_s\}$ gives a bounding
process, and the pair $(L_s, V_s)$ is a quite concise representation of~$D_s$. 

Recall that our construction~(\ref{drive}) of the Markov chain~$\vecX$ with
kernel~$K$ relies on the choice of a transition rule~$(\phi, \mu)$
satisfying~(\ref{kq}).

\begin{definition}
\label{mondef}
\emph{A transition rule~$\phi$ is said to be \emph{monotone\/} if each of the
mappings \linebreak $\phi(\cdot, u): \Xc \to \Xc$,\ \,$u \in \Uc$, is monotone
increasing, i.e.,\ if
$$
\mbox{$\phi(x, u) \leq \phi(y, u)$ for all $u \in \Uc$}
$$
whenever $x \leq y$.
}
\end{definition}

Suppose now that~$\phi$ is monotone.  Set $(L_0, V_0) := (\zh, \oh)$ and,
inductively,
$$
(L_s, V_s) := (\phi(L_{s - 1}, U_s), \phi(V_{s - 1}, U_s)),\ \ \ 1 \leq s \leq
t.
$$
One immediately verifies by induction that~(\ref{squeeze}) is satisfied.  Note
that~$(L_s, V_s)$ is
determined solely by~$\vecU$
(as is the coalescence detection
event $C = \{L_t = V_t\}$)
and is nothing more than $(Y_s(\zh), Y_s(\oh))$.  In plain language, since
monotonicity is preserved, when the chains~$\vecY(\zh)$ and~$\vecY(\oh)$ have
coalesced, so must have every~$\vecY(x)$.

\begin{remark}
\label{anti}
(a)~Lower and upper bounding processes can also be constructed when
Algorithm~\ref{nontech} is applied with a so-called ``anti-monotone'' transition
rule; we omit the details.  See H\"{a}ggstr\"{o}m and Nelander~\cite{HN},
Huber~\cite{Hexact}, Kendall~\cite{Kendall}, M{\o}ller~\cite{Moller}, M{\o}ller
and Schladitz~\cite{MS}, and Th\"{o}nnes~\cite{Thonnes} for further discussion in
various specialized settings.  There are at least two neat tricks associated with
anti-monotone rules.  The first is that, by altering the natural partial order
on~$\Xc$, such rules can be regarded, in certain bipartite-type settings, as
monotone rules, in which case the performance analysis in Section~5.3 of~\cite{FMMR} is 
available:\ consult Section~3 of~\cite{HN}, the paper~\cite{MS}, and Definition~5.1
in~\cite{Thonnes}.  The second is that the poset~$\Xc$ is allowed to be ``upwardly
unbounded'' and so need not have a~$\oh$:\ consult Section~2 of~\cite{Moller} and,
again, \cite{MS} and~\cite{Thonnes}.

(b)~Dealing with monotone rules on partially ordered state spaces without~$\oh$ is
problematic and requires the use of ``dominating processes.''
We comment that a dominating process provides a sort of {\em random\/} bounding
process and is useful when the state space is noncompact, but we shall not pursue
these ideas any further here. 
See
Kendall~\cite{Kendall} and Kendall and M{\o}ller~\cite{KM} in the context of CFTP;
we hope to discuss the use of dominating processes for our algorithm in future
work.
\end{remark}

\section{Our algorithm and~CTFP:\ comparison and connection}
\label{taleof2}

\subsection{Comparison}
\label{comparison}

How does our extension of Fill's algorithm, as given by Algorithm~\ref{nontech}
and discussed in detail in Section~\ref{app1}, compare to CFTP?  As we see it, our
algorithm has two main advantages and one main disadvantage.
\medskip

\emph{Advantages:\/}\ As discussed in Section~\ref{intro} and
Remark~\ref{genremark}(b) and in~\cite{Fill},
a primary
advantage of our algorithm
is interruptibility.
Given the close connection between the algorithms described in
Section~\ref{conn},
one may reasonably view the extra computational costs of our algorithm
(see ``{\em Disadvantage\/}'' below) as the costs of securing interruptibility.
A related second advantage concerns
memory allocation.  Suppose, for example, that our state space~$\Xc$ is finite and
that each time-step of Algorithm~\ref{nontech}, including the necessary imputation
(recall Section~\ref{imp}), can be carried out using a bounded amount of memory. 
Then, for fixed~$t$, our algorithm can be carried out using a fixed finite amount
of memory.  Unfortunately, it is rare in practice that the kernel~$K$ employed is
sufficiently well analyzed that one knows in advance a value of~$t$ (and a value
of the seed~$\XX_t$) giving a reasonably large probability~$\PP(\CC)$ of
acceptance.  Furthermore, the fixed amount of memory needed is in practice larger
than the typical amount of memory allocated dynamically in a run of CFTP.
Finally, we should note that Wilson~\cite{Wilson_read_once} has very recently
presented a version of CFTP which also can be carried out with a fixed finite amount
of memory, and which does not require an {\em a priori\/} estimate of the mixing time
of the chain.
\medskip

\emph{Disadvantage:\/}\ A major disadvantage of our Algorithm~\ref{nontech}
concerns computational complexity.  We refer the reader to~\cite{Fill}
and~\cite{MTFalgo} for a more detailed discussion in the setting of our
Section~\ref{mono} (and, more generally, the setting of stochastic monotonicity).
Briefly, if
no attention is paid to memory usage, our algorithm has running time competitive
with CFTP: cf.~Remark~\ref{genremark}(c), and also the discussion in Remark~9.3(e)
of~\cite{Fill} that the running time of our algorithm is, in a certain sense,
best possible in the stochastically monotone setting.  However, this analysis
assumes that running time is measured in Markov chain steps;  unfortunately,
time-reversed steps can sometimes take longer than do forward steps to execute
(e.g.,~\cite{MTFalgo}), and the imputation described in Section~\ref{imp} is
sometimes difficult to carry out.  Moreover, the memory usage for naive
implementation of our algorithm can be exorbitant; how to trade off speed for
reduction in storage needs is described in~\cite{Fill}.

\subsection{An alternative to Algorithm~1.1}
\label{altalgsub}

Thus far we have been somewhat sketchy about the choice(s) of~$t$ in
Algorithm~\ref{nontech}.  As discussed in Section~\ref{intro}, one possibility is
to run the repetitions of the basic routine independently, doubling~$t$ at each
stage.  However, another possibility is to continue back in time, reusing the
already imputed values~$\UU_s$ and checking again for coalescence.  (There is an
oblique reference to this alternative in Remark~9.3 of Fill~\cite{Fill}.)  This
idea leads to the following algorithm.

\begin{algorithm}
\label{altalg}
Choose an initial state $\XX_0 \sim \pih$, where $\pih$ is absolutely
continuous with respect to~$\pi$.  Run the time-reversed chain~$\Kt$, obtaining
$\XX_0, \XX_{-1}, \ldots$ in succession.  Conditionally given $(\XX_0,
\XX_{-1}, \ldots) = (x_0, x_{-1}, \ldots)$, generate independent random variables
$\UU_0, \UU_{-1}, \ldots$ with marginals
\begin{equation}
\label{altimp}
\PP(\UU_s \in du) = P(U \in du\,|\,\phi(x_{s - 1}, U) = x_s),\ \ \ s = 0, -1,
\ldots,
\end{equation}
where, on the right, $\Lc_P(U) = \mu$ is given by~(\ref{kq}).  For $t = 0, 1,
\ldots$ and $x \in \Xc$, set $\YY^{(-t)}_{-t}(x) := x$ and, inductively,
$$
\YY^{(-t)}_s(x) := \phi(\YY^{(-t)}_{s - 1}(x), \UU_s),\ \ \ -t + 1 \leq s \leq 0.
$$
If $\TT < \infty$ is the smallest~$t$ such that
\begin{equation}
\label{altcoalescence}
\mbox{$\YY^{(-t)}_0(x)$, $x \in \Xc$, agree\ \ \ \ \ \ \ ($ = \XX_0$),}
\end{equation}
then the algorithm succeeds and reports~$\WW := \XX_{-\TT}$ as an observation
from~$\pi$.  Otherwise, the algorithm fails.
\end{algorithm}

\begin{remark}
\label{altremark}
(a)~We need only generate $\XX_0, \XX_{-1}, \ldots, \XX_{-t}$ and then
impute $\UU_0$, $\UU_{-1}, \ldots, \UU_{-t + 1}$
using~(\ref{altimp}) in order to check
whether or not~(\ref{altcoalescence}) holds.  Thus if $\TT < \infty$, then the
algorithm terminates in finite time.

(b)~We omit the detailed description \`{a} la Section~\ref{app1}.
But the key in setting up the first probability space is \emph{first\/} to choose $W \sim
\pi$ and $U_0, U_{-1}, \ldots$ all mutually independent and \emph{then\/}, having
determined the backwards coalescence time~$T$ from $U_0, U_{-1}, \ldots$, to set $X_{-T}
:= W$.

(c)~We may relax the condition that~$\TT$ be the
\emph{smallest\/}~$t$ satisfying~(\ref{altcoalescence}), via the use of coalescence
detection events as in Section~\ref{app2}.  In particular, to save considerably
on computational effort, we may let~$\TT'$ be the smallest~$t$ which is a power
of~$2$ such that~(\ref{altcoalescence}) holds and report $\XX_{-\TT'}$ instead.
%

(d)~Algorithm~\ref{altalg}, and likewise its variant in remark~(c), is
interruptible:\ \ $\TT$ and~$\WW$ are conditionally independent given success.

(e)~Uniform ergodicity of~$K$ is necessary (and, in a weak sense, sufficient) for
almost sure success of Algorithm~\ref{altalg}; cf.\ Remark~\ref{genremark}(a).
\end{remark}

\subsection{Connection with CFTP}
\label{conn}

There is a strong and simple connection between CFTP and our
Algorithm~\ref{altalg}.  Indeed, suppose we carry out the usual CFTP algorithm to
sample from~$\pi$, using kernel~$K$, transition rule~$\phi$, and driving
variables~$\vecU = (U_0, U_{-1}, \ldots)$.  Let~$T$ denote the backwards
coalescence time and let $X_0 \sim \pi$ denote the terminal state output by CFTP. 
Let $W \sim \pi$ independent of~$\vecU$, and follow the trajectory from $X_{-T} :=
W$ to $X_0$; call this trajectory $\vecX = (X_{-T}, \ldots, X_0)$.  Since~$X_0$ is
determined solely by~$\vecU$, the random variables~$W$ and~$X_0$ are independent.

When $\pih = \pi$ in Algorithm~\ref{altalg}, the algorithm simply constructs the
same probability space as for CFTP, but with the ingredients generated in a
different chronological order:\ first $X_0, X_{-1}, \ldots$; then~$\vecU$ (which
determines~$T$); then $W := X_{-T}$.  Again $X_0 \sim \pi$ and $W \sim \pi$ are
independent.

\begin{remark}
\label{connremark}
(a)~Because of this statistical independence, it does not matter in
Algorithm~\ref{altalg} that we actually use $\XX_0 \sim \pih \neq \pi$.

(b)~The fact (1)~that $W$, unlike~$X_0$, is independent of~$\vecU$, together with
(2)~that $T$ depends solely on~$\vecU$, explains why our algorithm is interruptible
and CFTP is not.

(c)~In a single run of CFTP, the user would of course be unable to choose $W
\sim \pi$ as above, just as in a single run of Algorithm~\ref{altalg} we do not
actually choose $X_0 \sim \pi$.  So one might regard our described connection
between the two algorithms as a bit metaphorical.  But see Section~7.2 of~\cite{FMMR}.
\end{remark}

\reversemarginpar

\end{document}